# An Upwind Finite Difference Method to Singularly Perturbed Convection Diffusion Problems on a Shishkin Mesh.

Daniel T. Gregory

**Abstract:**

This paper introduces a numerical approach to solve singularly perturbed convection diffusion boundary value problems for second-order ordinary differential equations that feature a small positive parameter ε multiplying the highest derivative. We specifically examine Dirichlet boundary conditions. To solve this differential equation, we propose an upwind finite difference method and incorporate the Shishkin mesh scheme to capture the solution near boundary layers. Our solver is both direct and of high accuracy, with computation time that scales linearly with the number of grid points. MATLAB code of the numerical recipe is made publicly available. We present numerical results to validate the theoretical results and assess the accuracy of our method. The tables and graphs included in this paper demonstrate the numerical outcomes, which indicate that our proposed method offers a highly accurate approximation of the exact solution.

---

[1] Department of Mathematics, University of Kentucky, Lexington, KY, 40508, e-mail: dannytgregory@uky.edu

## 1. Introduction:

In this paper, we consider the following stationary state second order singularly perturbed differential equation with Dirichlet boundary conditions.

$$L_\varepsilon \equiv -\varepsilon u''(x) + a(x)u'(x) = f(x) \qquad (1.1)$$

$$u(0) = 0 \qquad u(1) = 0$$

where, $u'' = \frac{d^2u}{dx^2}$ is the second derivative of $u$ with respect to the independent variable $x$, and $f(x)$ is called a load/source function. The solution $u$ of the equation is an unknown scalar function. There is no initial condition because the equation does not depend on time hence it becomes a boundary value problem. The function $f(x)$, and $a(x)$, are all smooth and satisfy $a(x) \geq \alpha > 0$. The parameter $\varepsilon$ is assumed to be a small positive value, such that $0 < \varepsilon \leq 1$. While it may be easy to analytically solve this problem in some cases, finding the solution $u$ with analytical techniques can be difficult or even impossible for a general function $f$ and/or with complicated boundary conditions. The well-posedness of problem (1.1) has been discussed in more detail in [1].

The Convection-Diffusion formula finds applications in various real-world scenarios such as predicting groundwater pollution [7], flows in chemical reactors [8], convective heat transport problems [9], and simulation of oil extraction from underground reservoirs[10]. For example the term "convection" refers to the upward movement of hotter and denser substances due to their heat, while colder and less dense substances tend to sink. This motion of substances leads to the transfer of heat. Meanwhile, diffusion occurs simultaneously with convection and involves the spreading of particles from areas of high concentration to regions of low concentration.

Several methods can be employed to solve these problems, including , finite volume method, finite element method and simple upwind finite difference method [4], [5], [ 6]. Unfortunately, a simple up-wind scheme can not capture the features of the boundary layer exactly because of the pollution effects. The problem being studied in this paper involves rapidly changing solutions in very thin regions near the boundary. Traditional numerical methods often fail to accurately capture these changes, which can result in errors across the entire domain. To address this issue, various methods such as Bakhavalov [11] and Gartland meshes [12] have been developed. In this study, we analyze a standard upwind finite difference method combined with the Shishkin mesh, which is a type of local refinement strategy introduced by a Russian mathematician Grigorii Ivanovich Shishkin in 1988 [13], [3].

**Remark 1.1:** For the sake of brevity, this paper primarily focuses on analyzing a one-dimensional problem, assuming a constant coefficient for the function $a(x)$ and a smooth function $f(x)$. In the subsequent discussion, the symbol $C$ represents a generic positive constant, which may vary across different formulas but remains unaffected by the mesh or the small positive parameter $\varepsilon$. A subscripted C (e.g., $C_1$) is also a constant that is independent of $\varepsilon$ and of any mesh used, but it takes one fixed value.

The subsequent sections of the article are structured as follows: Section 2 serves as an introduction to the finite difference method applied to equation (1.1), accompanied by the introduction of the layer-adapted Shishkin mesh. Moving forward, Section 3 presents error estimate theorems pertaining to both uniform mesh upwind finite difference schemes and upwind finite difference schemes employing the Shishkin mesh. In Section 4, we provide a comprehensive presentation of numerical results that serve to validate the theoretical findings we have discussed thus far. Finally, in Section 5, we conclude our results while also considering potential avenues for future research and exploration.

## 2. Finite Difference Method:

The Finite Difference Method (FDM) is a numerical technique used to approximate solutions to differential equations. It is particularly useful when analytical solutions are difficult or impossible to obtain. The derivatives in the original differential equation are approximated using finite difference formulas. By substituting these approximations into the differential equation, the problem is transformed into a system of algebraic equations. There are different types of finite difference approximations, such as forward difference, backward difference, and central difference. The choice of the specific formula depends on the nature of the problem and the desired accuracy. Moreover, finite difference formulation offers a more direct and intuitive approach to the numerical solution of differential equations than other formulations.

### 2.1 Finite Difference Discretization

Consider the purely singularly perturbed convection-diffusion problem
$$L_\varepsilon \equiv -\varepsilon u''(x) + a(x)u'(x) = f(x) \quad \text{for} \quad 0 < x < 1 \quad (2.1)$$
$$u(0) = 0 \quad u(1) = 0.$$

where, $a(x) > \alpha > 0$, and $0 < \varepsilon \leq 1$. Assume that $a(x)$ and $f(x)$ lie in $c^{\infty}[0, 1]$. Let N be a positive integer. Partition $[0, 1]$ by the equidistant mesh $x_i = ih$ for $i = 0, 1,..., N$, where $h := \frac{1}{N}$. Standard discretizations of differential equations use a central difference approximation. That is from the Taylor series expansion we can obtain the following central difference approximations for derivatives.

$$u'(x) \approx \frac{u_{i+1} - u_{i-1}}{2h} \quad \text{and} \quad u''(x) \approx \frac{u_{i-1} - 2u_i + u_{i+1}}{h^2}$$

for $i = 1, , ... , N - 1$.

For the sake of simplicity, let's consider the function $a(x) = 1$. In this case, we can express the discretized version of equation (2.1) as a matrix equation, as follows:

$$(-\varepsilon A + B)U = F \qquad (2.2)$$

where, A represent the stiffness matrix corresponding to $u''(x)$, B represent the convection matrix corresponding to $u'(x)$, F represent the source term corresponding to f, and U represent the discrete solution of equation (2.1) as follows:

$$A = \begin{bmatrix} -2 & 1 & & & \\ 1 & -2 & 1 & & \\ & \ddots & \ddots & \ddots & \\ & & 1 & -2 & 1 \\ & & & 1 & -2 \end{bmatrix} \qquad U = \begin{bmatrix} u_1 \\ u_2 \\ \vdots \\ u_{n-1} \end{bmatrix}$$

$$B = \frac{1}{2}\begin{bmatrix} 0 & 1 & & & \\ -1 & 0 & 1 & & \\ & \ddots & \ddots & \ddots & \\ & & -1 & 0 & 1 \\ & & & -1 & 0 \end{bmatrix} \qquad F = \begin{bmatrix} f_1 \\ f_2 \\ \vdots \\ f_{n-1} \end{bmatrix}$$

**Remark 2.1:** For singularly perturbed problem there exist $\hat{x} \in [0, 1]$ such that,

$$\lim_{\varepsilon \to 0} \lim_{x \to \hat{x}} u(\hat{x}) \neq \lim_{x \to \hat{x}} \lim_{\varepsilon \to 0} u(\hat{x}) \tag{2.3}$$

One can easily check that for $f(x) = x$, and $a(x) = 1$ the equation (2.1) has the exact solution

$$u(x) = x\left(\frac{x}{2} + \varepsilon\right) - \left(\frac{1}{2} + \varepsilon\right)\left(\frac{e^{(x-1)/\varepsilon} - e^{-1/\varepsilon}}{1 - e^{-1/\varepsilon}}\right) \tag{2.4}$$

Moreover, at $x = 1$ the solution $u(x)$ satisfies the condition (2.3) and thus the boundary value problem (2.1) has a boundary layer at $x = 1$. This is a narrow region where $u$ is bounded independently of $\varepsilon$ but its derivatives blow up as $\varepsilon \to 0$.

### 2.2 Layer Adapted Shishkin Mesh

In this section we introduce the Shishkin mesh named after the Russian mathematician G. I. Shishkin [14], [13], [3]. The purpose of using a Shishkin mesh is to accurately capture the behavior of a solution $u(x)$ near regions where the solution exhibits sharp changes or boundary layers. Typically, a Shishkin mesh is constructed by dividing the domain into different subdomains and then applying a different mesh spacing in each subdomain. The mesh spacing is chosen such that it becomes finer as the solution approaches regions of interest, such as boundary layers or areas of rapid variation. This adaptive meshing strategy helps to ensure that the numerical approximation captures the important features of the solution more effectively.

Let, $X_s^N : 0 = x_0 < x_1 < x_2 < \ldots < x_{n-1} < x_n = 1$ be the Shishkin mesh constructed by dividing the interval $[0, 1]$ into two subintervals for an even positive integer $N$.

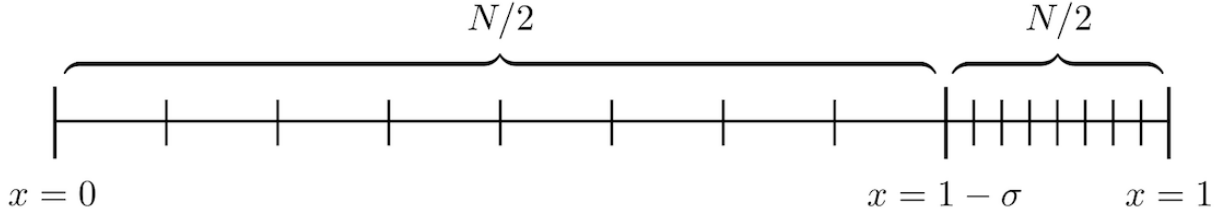

Figure 1: A Shishkin mesh refinement

where,

$$\sigma = \min\left\{\frac{1}{2}, \left(\frac{2}{\alpha}\right)\varepsilon \ln N\right\}$$

which depends on ε and *N*. In our exposition we shall assume that σ = (2/α)εlnN, as σ = 1/2 occurs only when N is exponentially large relative to ε, which is rare in practice. Then the mesh transition point, which separates the fine and coarse portions of the Shishkin mesh, is defined to be 1 − σ; typically it lies close to 1. For an even integer N, divide each of [0, 1 − σ] and [1 − σ, 1] by an equidistant mesh with N/2 subintervals.

Graded meshes, where the mesh width gets finer and finer as one moves closer and closer to x = 1, have been advocated by several authors; see [1] for references. But convergence analyses on graded meshes can be very delicate, so we shall concentrate here on a simpler piecewise-equidistant mesh that is originally due to Shishkin and which has been used by many researchers. See [15, 16] for references to singularly perturbed problems where the Shishkin mesh has been used.

### 3. Error Estimates

Solutions to singularly perturbed convection diffusion equations can be written as a well-behaved term and a layer term. Such decompositions of *u* are helpful when constructing accurate numerical methods and are often needed in the error analysis.

**Definition 3.1:** Let *u* be a function defined on a domain $D \subset [0, 1]$. We define $L_\infty$ norm or maximum norm by $||u||_\infty = \max|u(x)|$ for $x \in \overline{D}$.

In the rest of this paper we assume that each norm of *u* is well-defined and finite.

**Theorem 3.1:** Let u be the solution to equation (2.1). Let $p$ be a nonnegative integer. Then there is a splitting $u = u_S + u_E$ such that, for $0 \leq j \leq p$, the inequalities

$$|u_S^j| \leq C \quad \text{and} \quad |u_E^j|_\infty \leq C\varepsilon^{-j} e^{-\alpha(1-x)/\varepsilon} \quad \text{for} \quad 0 \leq x \leq 1$$

hold for some constant $C = C(p)$, where $u_S$ is the smooth part of the solution and $u_E$ is the layer part of the solution.

*proof:* The proof can be found from the reference [2].

### 3.1 Upwind Scheme

In computational physics, the term upwind scheme typically refers to a class of numerical discretization methods for solving differential equations, in which so-called upstream variables are used to calculate the derivatives in a flow field. In this work, upwinding means taking a one-sided difference on the side away from the layer. i.e.,

$$u'(x) \approx \frac{u_i - u_{i-1}}{h}$$

In this way we can avoid unnatural oscillations from the computed solution. However, the layers in the computed solution are excessively smeared, meaning they are not as steep as they should be. For a rich study of upwinding please refer [6].

**Theorem 3.2** Let $\{u_i^N\}_{i=0}^N$ be the solution to (2.1) computed using simple upwinding on a equidistant mesh with $N$ subintervals. Suppose that $h \gg \varepsilon$. Then there exist a constant C such that

$$|u_i - u_i^N| \leq C[h + exp\left(\frac{-\alpha(1-x_i)}{\alpha h + 2\varepsilon}\right)] \text{ for } i = 0, \ldots, N$$

if $x_i$ is bounded away from $x = 1$, then $|u_i - u_i^N| \leq C[h + exp\left(\frac{-\alpha(1-x_i)}{\alpha h + 2\varepsilon}\right)] \leq ch$

That is, the upwind scheme yields an $o(h)$-accurate solution away from, But at interior mesh points that lie close to or inside the layer the scheme is only o(1)-accurate.

*Proof:* The proof can be found from [2].

### 3.2 Estimate on Shishkin mesh

To analyze the convergence of the solution we split the discrete solution $u_i^N = u_{S,i}^N + u_{E,i}^N$ and then by the triangle inequality we get,

$$|u_i - u_i^N| = |(u_S + u_E)_i - (u_{S,i}^N + u_{E,i}^N)| \leq |u_{S,i} - u_{S,i}^N| + |u_{E,i} - u_{E,i}^N| \quad (3.1)$$

for all $i = 0, \ldots, N$. Now we shall bound each term of the equation (3.1) above separately.

**Lemma 3.1** There exist a constant $C$ such that

$$|u_{S,i} - u_{S,i}^N| \leq C N^{-1}. \quad i = 0, \ldots, N.$$

*Proof.* The proof can be found from [2].

**Lemma 3.2** There exist constants $C_0$ and $C_1$ such that

$$|u_{E,i} - u_{E,i}^N| \leq C_0 N^{-1}. \quad i = 0, \ldots, N/2$$

$$|u_{E,i} - u_{E,i}^N| \leq C_0 N^{-1} lnN. \quad i = N/2 + 1, \ldots, N.$$

*Proof.* The proof can be found from [2].

By Lemmas 3.1 and 3.2 the final convergence result under shishkin mesh can now be stated as follows:

**Theorem 3.3:** There exist a constant $C$ such that the discrete solution $u_i^N$ of (2.2) under the finite difference method satisfies

$$|u_i - u_i^N| \leq C N^{-1} lnN. \quad i = 0, \ldots, N$$

*Proof.* Combine the equation (3.1) with Lemma 3.1 and lemma 3.2.

## 4. Numerical Results

In this section we present a few numerical experiments to illustrate the computational method discussed in this paper. The numerical experiments are performed on a laptop computer with MATLAB R2022a in MacBook Air with M1 chip. We use the following formula for the order of convergence.

$$\text{Order of Convergent} = \frac{\log\left((u_{i+1} - U_{i+1})/(u_i - U_{i+1})\right)}{\log(2)} \qquad (4.1)$$

where $u_{i+1}$ is the numerical solution in the $(i+1)^{th}$ mesh iteration and $U_{i+1}$ is the analytical solution computed at the same mesh level. Here we consider the following model problem to validate theoretical results.

$$\begin{aligned} -\varepsilon u'' + u' &= x \\ u(0) &= 0 \\ u(1) &= 0. \end{aligned} \qquad (4.2)$$

The exact solution is $u(x) = x\left(\frac{x}{2} + \varepsilon\right) - \left(\frac{1}{2} + \varepsilon\right)\left(\frac{e^{(x-1)/\varepsilon} - e^{-1/\varepsilon}}{1 - e^{-1/\varepsilon}}\right) \qquad (4.3)$

**Example 1:** In this example we solve equation (4.2) to obtain the numerical solution and compare it with the exact solution in the equation (4.3). We compare maximum norm errors between numerical solution $\left(U^N\right)$ and exact solution $\left(u^N\right)$ solutions computed on same mesh levels obtained through sequences of meshes refinements to check the accuracy of the method where $N$ is the number of nodal points on the domain.

**Test case 1:** We first consider solving equation (4.2) by FDM using central difference scheme with $\varepsilon = 1$. The resulting maximum norm errors for different mesh sizes (N) are presented in Table 1. The table 1 and figure 2 conclude that the numerical solution effectively approximates the exact solution with a high accuracy. It is worth noting that there is no need to employ an upwind scheme or Shishkin mesh in this case, as the chosen value of $\varepsilon = 1$ does not induce the presence of a boundary layer

Table 1: Maximum norm error for $\varepsilon = 1$ under uniform mesh refinements

| N | 256 | 512 | 1024 | 2048 | 4096 | 8192 | 16384 |
|---|---|---|---|---|---|---|---|
| $\|\|u^N - U^N\|\|_\infty$ | 2.28e-07 | 5.73e-08 | 1.43e-08 | 3.59e-09 | 8.99e-10 | 2.70e-10 | 5.50e-11 |

Figure 2: Solution plots for $\varepsilon = 1$ under a uniform mesh with central difference scheme

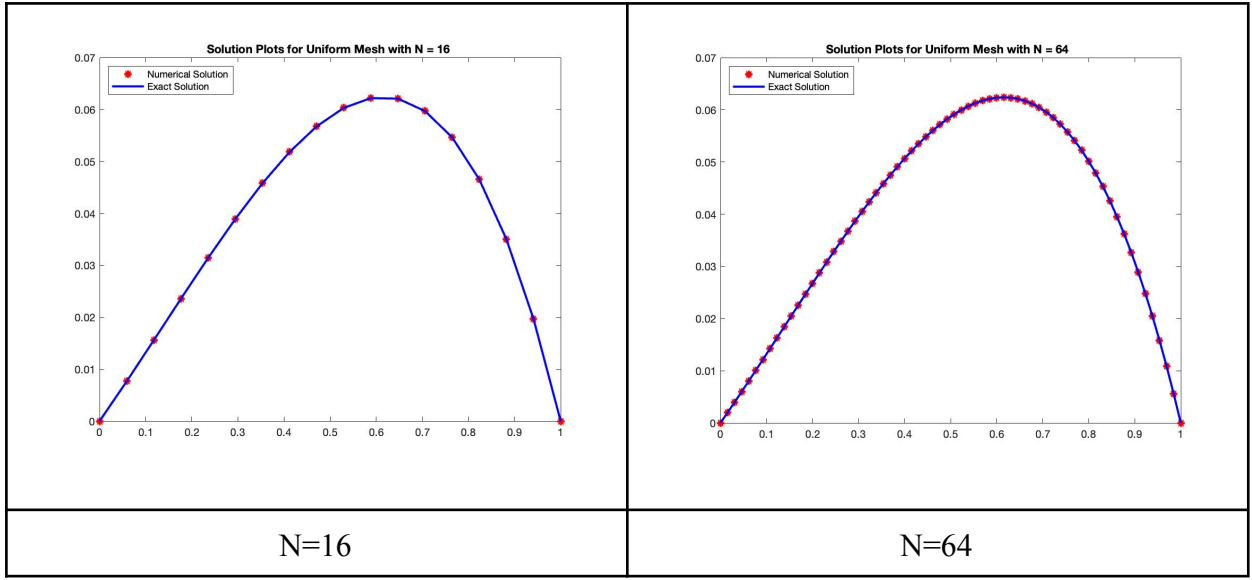

| N=16 | N=64 |

**Test case 2:** Secondly, we solve equation (4.2) for $\varepsilon = 10^{-2}, 10^{-4}, 10^{-6}, 10^{-8}$ through a finite difference scheme with central differences under uniformly refined meshes. We then compare the obtained numerical solution with the exact solution and record the maximum norm error in the table 2. Table 1 and Figure 2 conclude that the numerical solution does not approximate the exact solution accurately under a uniform mesh when the problem shows high singularly perturbed features.

Table 2: The maximum norm error $\|\|u^N - U^N\|\|_\infty$ under uniform meshes with central difference

| $\varepsilon \setminus N$ | 256 | 512 | 1024 | 2048 | 4096 | 8192 | 16384 |
|---|---|---|---|---|---|---|---|
| $\varepsilon = 10^{-2}$ | 0.0024 | 5.96e-04 | 1.48e-04 | 3.72e-05 | 9.31e-06 | 2.32e-06 | 5.82e-07 |
| $\varepsilon = 10^{-4}$ | 0.4512 | 0.4070 | 0.3300 | 0.2132 | 0.0932 | 0.0265 | 0.0058 |
| $\varepsilon = 10^{-6}$ | 0.4994 | 0.4985 | 0.4960 | 0.4959 | 0.4919 | 0.4839 | 0.4683 |

| $\varepsilon = 10^{-8}$ | 0.5000 | 0.5000 | 0.5000 | 0.5000 | 0.4999 | 0.4997 | 0.4997 |

Figure 3a: Solution plots under uniform meshes under central difference scheme

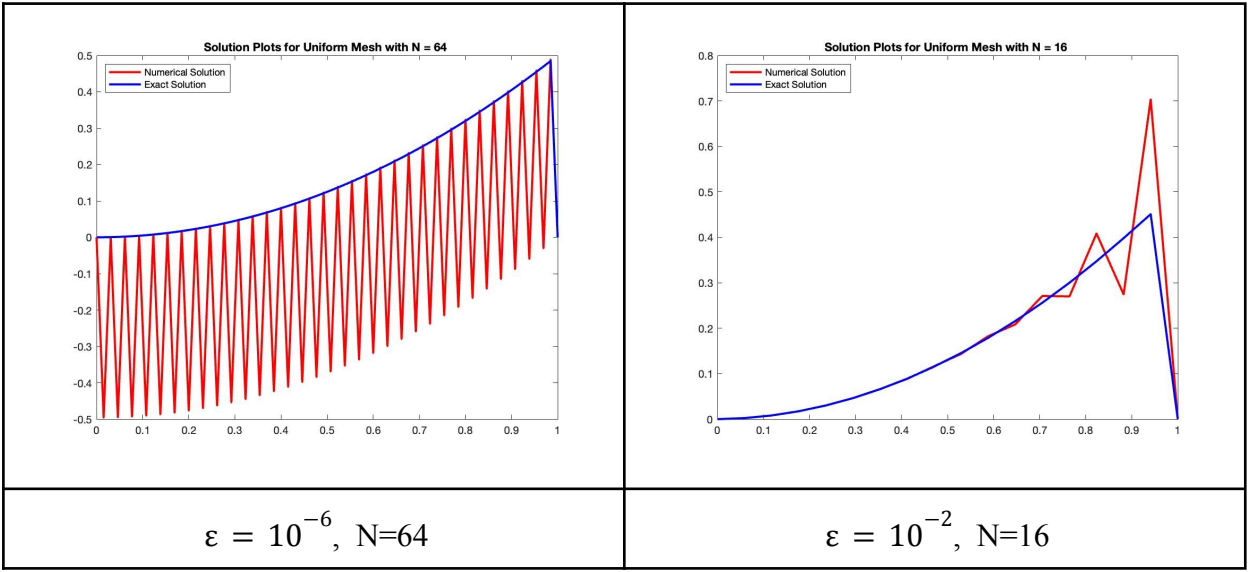

| $\varepsilon = 10^{-6}$, N=64 | $\varepsilon = 10^{-2}$, N=16 |

Figure 3b: Solution plots under uniform meshes under central difference scheme

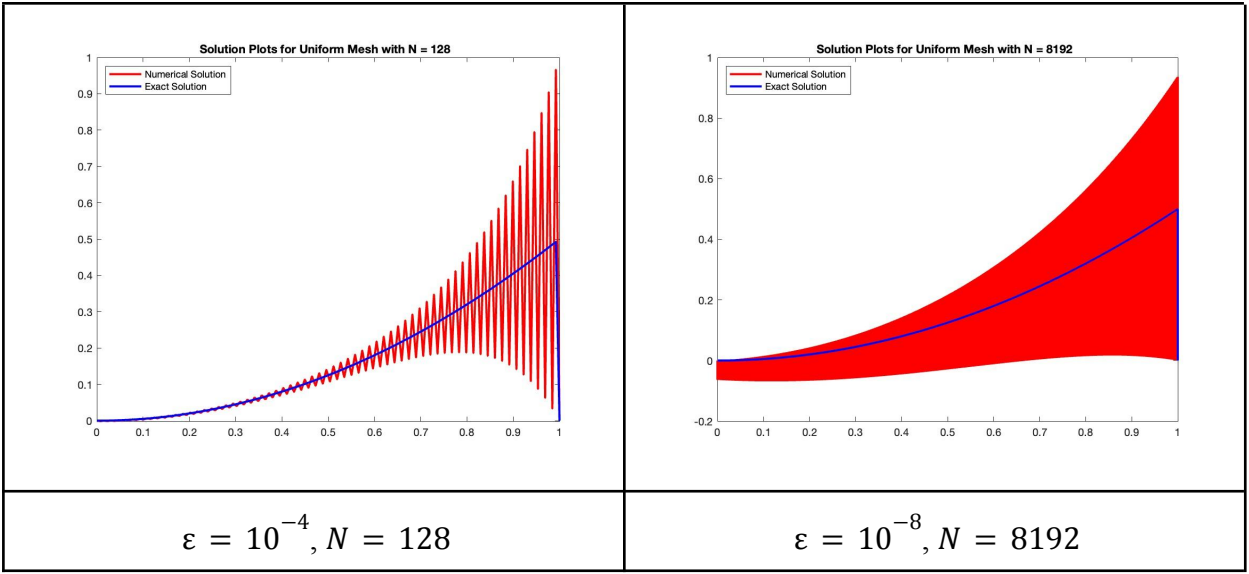

| $\varepsilon = 10^{-4}, N = 128$ | $\varepsilon = 10^{-8}, N = 8192$ |

**Test case 3:** Finally, to overcome the issue that arises in **Test case 2** we repeat **Test case 2** with an upwinding finite difference scheme under Shishkin mesh. Table 3 and figure 4 conclude that the numerical solution approximate the exact solution accurately in the presence of high singularly perturbed features.

Table 3: The maximum norm error $||u^N - U^N||_\infty$ under Shishkin meshes with upwind scheme

| ε \ N | 256 | 512 | 1024 | 2048 | 4096 | 8192 | 16384 |
|---|---|---|---|---|---|---|---|
| $\varepsilon = 10^{-2}$ | 0.0118 | 0.0091 | 0.0075 | 0.0067 | 0.0063 | 0.0060 | 0.0057 |
| $\varepsilon = 10^{-4}$ | 0.0063 | 0.0038 | 0.0022 | 0.0013 | 7.38e-04 | 4.35e-04 | 2.66e-04 |
| $\varepsilon = 10^{-6}$ | 0.0063 | 0.0037 | 0.0021 | 0.0012 | 6.78e-04 | 3.74e-04 | 2.04e-04 |
| $\varepsilon = 10^{-8}$ | 0.0063 | 0.0037 | 0.0021 | 0.0012 | 6.78e-04 | 3.74e-04 | 2.04e-04 |

Figure 4a: Solution plots under Shishkin mesh with upwind scheme

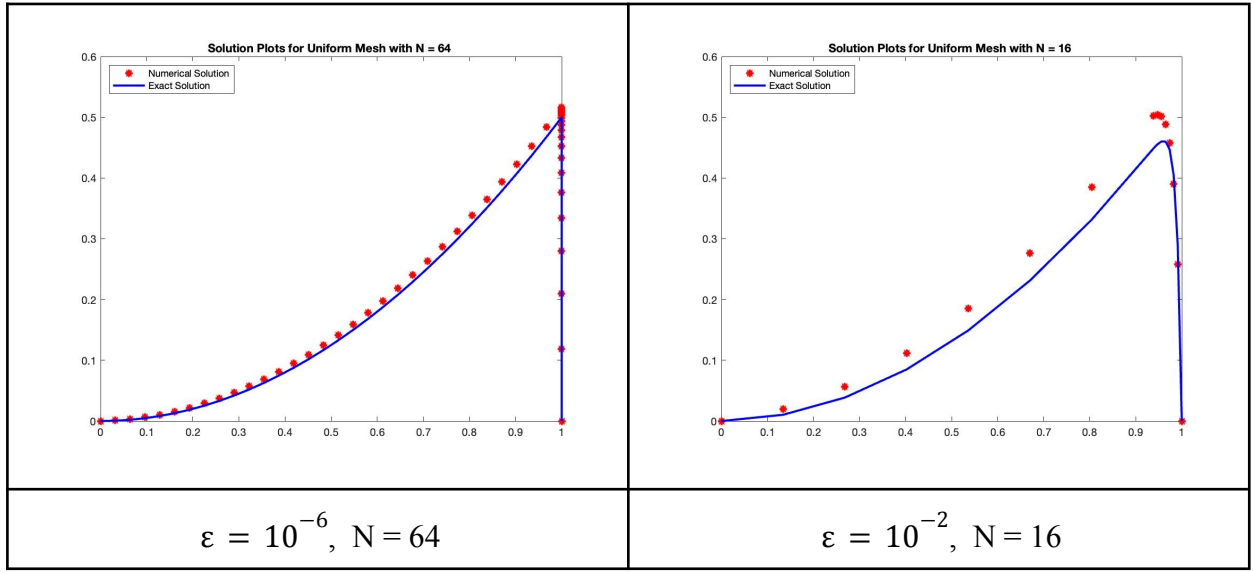

| $\varepsilon = 10^{-6}$, N = 64 | $\varepsilon = 10^{-2}$, N = 16 |

Figure 4b: Solution plots under Shishkin mesh with upwind scheme

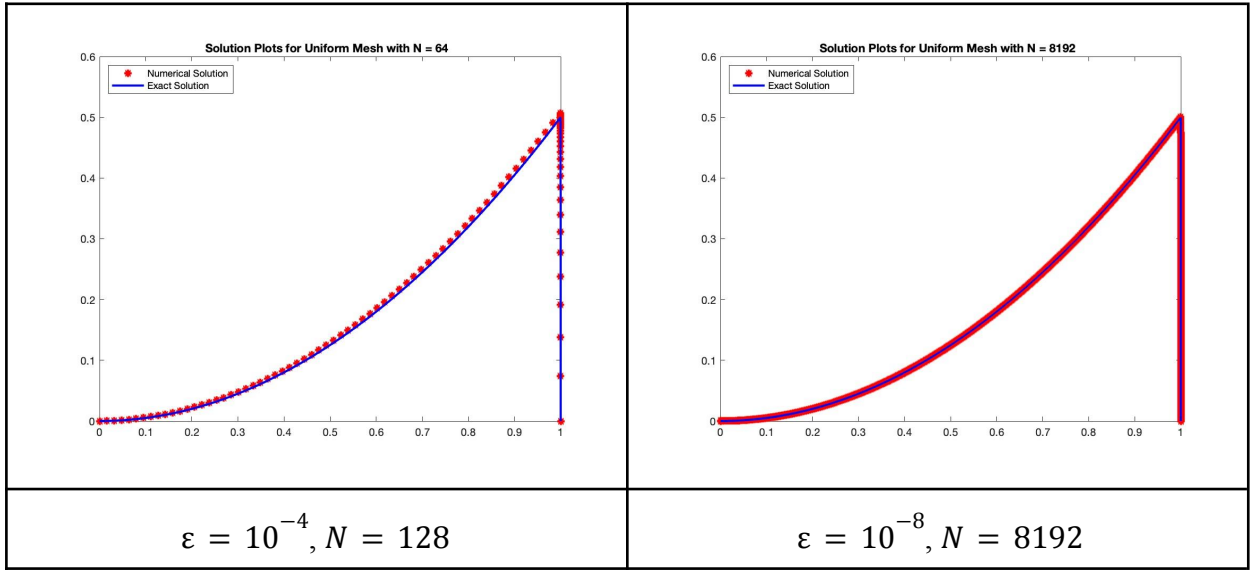

| $\varepsilon = 10^{-4}, N = 128$ | $\varepsilon = 10^{-8}, N = 8192$ |

**Example 2:** In table 3, we record convergence rates of the numerical solution obtained through a sequence of meshes used to validate results in the theorem 3.3. For $\varepsilon = 10^{-2}, 10^{-4}, 10^{-6},$ and $10^{-8}$ an upwinding Shishkin mesh is used to attack the boundary layer near $x = 1$. Convergence rates are closely matched with the results in theorem 3.3. In particular, for $\varepsilon = 1$, a uniformly refined mesh is used since the problem has no boundary layer. It is shown that when $\varepsilon = 1$ convergence rates are close to 2 which perfectly match with the standard theoretical results.

Table 4: The maximum norm error $||u^N - U^N||_\infty$ under Shishkin meshes with upwind scheme

| N | 256 | 512 | 1024 | 2048 | 4096 | 8196 | 16384 |
|---|---|---|---|---|---|---|---|
| $\varepsilon = 1$ | 1.9888 | 1.9944 | 1.9972 | 1.9986 | 1.9991 | 1.7359 | 2.2941 |
| $\varepsilon = 10^{-2}$ | 0.4443 | 0.3771 | 0.2736 | 0.1665 | 0.0922 | 0.0650 | 0.0668 |
| $\varepsilon = 10^{-4}$ | 0.6882 | 0.7417 | 0.7752 | 0.7906 | 0.7886 | 0.7634 | 0.7053 |
| $\varepsilon = 10^{-6}$ | 0.6932 | 0.7505 | 0.7912 | 0.8201 | 0.8414 | 0.8576 | 0.8700 |

| $\varepsilon = 10^{-8}$ | 0.6933 | 0.7506 | 0.7913 | 0.8204 | 0.8420 | 0.8587 | 0.8719 |

**Example 3:** In this example we compute CPU time required by the upwind finite difference method (FDM) to observe the efficiency of the upwind finite difference method for a sequence of meshes. In particular we compare CPU time for $\varepsilon = 1, 10^{-8}$. This concludes that the upwind finite difference method works fast on 1D problems due to its simple structure.

Table 5: CPU time in seconds

| N | 128 | 256 | 512 | 1024 | 2048 | 4096 | 8192 | 16384 |
|---|---|---|---|---|---|---|---|---|
| $\varepsilon = 1$ | 0.0500 | 0.0600 | 0.0800 | 0.1300 | 0.5100 | 2.1200 | 9.4900 | 43.4300 |
| $\varepsilon = 10^{-8}$ | 0.0600 | 0.1000 | 0.1300 | 0.1900 | 0.4600 | 1.2800 | 8.1200 | 41.0800 |

## 5. Conclusion and Future Work

In this paper we present an accurate and fast numerical method to solve one dimensional second order singularly perturbed convection diffusion equation with Dirithlet boundary condition. Due to the boundary layer that occurs at the right end point ($x = 1$), the standard central difference method with uniformly refined meshes does not approximate the exact solution accurately. To attack the oscillations near at $x = 1$, we introduced an upwind finite difference method under Shishkin mesh refinements. In this way, numerical solutions approximate the exact solutions accurately and efficiently as concluded in section 4.

This work can be extended to higher order one dimensional singularly perturbed differential equations with different boundary conditions like Dirichlet, Neumann or Robin boundary conditions [17]. The analysis could also be applied to higher dimensions and/or variable coefficients, although this may present some challenges. Singularly perturbed features and the method discussed in this paper can be applied to the work done by [18], [19], and [20]. Additionally, the problem could be expanded with singular source term $f$ with different numerical techniques like the finite element method.

**Acknowledgments:**

The author wishes to express gratitude to his project advisor, Dr. Charuka D. Wickramasinghe. Dr. Wickramasinghe played a crucial role in guiding the author through his initial undergraduate research project in mathematics and consistently provided support and encouragement.

**Appendix: MATLAB code**

We wrote the following Matlab code to illustrate the implementation of layer adapted Shishkin mesh and a uniform mesh.

This code implements the 1D Finite Difference Method for solving the singularly perturbed convection diffusion problem -epsilon*u" + u' = f with u(0) = 0 and u(1) = 0. It considers the problem with Dirichlet boundary conditions and provides solutions for both uniform mesh and Shishkin meshes. Finally, the code computes the maximum norm error by comparing the obtained solution with the exact solution. The exact solution is calculated based on the given problem equation. The maximum norm error is displayed as the output.

```
close;clc; format short
tStart = cputime;
epsilon=10^-6;  % Singularly perturbed parameter
p=[64 128 256 512 1024 2048 4096]; % Generate different meshes
for j=1:length(p) % a for loop to run through all the meshes are record errors and rates
%% Generate a uniform mesh
% n=p(j); % Number of interior points
% h=1/(n+1); % Step size
% y=0:h:1; % The interval (the domain)
% x=y';
% %% Generate a Shishkin mesh
```

```matlab
beta=0.9; sigma=2;
N=p(j); % number of nodal point in the domain (an even number)
tau = min(1/2, sigma*epsilon*log(N)/beta);
  x = unique([linspace(0,1-tau, N/2), ...
    linspace(1-tau,1, N/2)])';
n = length(x)-2;
h=zeros(n,1);
for i = 1:n-1
h(i) = x(i+1) - x(i); % step size in the Shishkin mesh
end
h(end)=h(end-1);
k=h.*h;
```

## Build matrices to get finite difference solution

```matlab
A =(2*diag(ones(1,n))) - 1*diag(ones(1,n-1),1) - 1*diag(ones(1,n-1),-1)); % Stiffness matrix
%C =0.5*h.*(diag(zeros(1,n)) + 1*diag(ones(1,n-1),1) - 1*diag(ones(1,n-1),-1));%Convection with central difference
C =h.*(1*diag(ones(1,n)) - 0*diag(ones(1,n-1),1) -1*diag(ones(1,n-1),-1));%Convection with upwind
f=k.*x(2:end-1); % Source term for f(x)=x
u_int=(epsilon*A+C)\f; % Finite difference approximation for interior points
U=zeros(n+2,1); % Sparse matrix to save all the solutions
U(1)=0; % Left boundary condition
U(end)=0; % Right boundary condition
U= [U(1); u_int; U(end)] ; % Complete numerical approximations
tEnd = cputime - tStart
```

## Exact solution for f(x)=x

```matlab
v1=x.*(0.5*x + epsilon);v2=(0.5+ epsilon);
v3=exp((x-1)/epsilon) - exp(-1/epsilon); v4=1 - exp(-1/epsilon);
u_exact= v1-v2.*(v3/v4);
```

## Record maximum norm errors

```matlab
evec=U-u_exact; maxerr=norm(evec, inf);
```

```
if j>2
rate= abs(log(errold/maxerr)/log(2));
end
errold=maxerr;
%% Plot solutions
figure(j)
plot(x,U,'r',x,u_exact,'b','LineWidth',2);
legend('Numerical Solution', 'Exact Solution','Location','NorthWest')
title('Solution Plots')
end
```